# LOCAL TAIL BOUNDS FOR FUNCTIONS OF INDEPENDENT RANDOM VARIABLES

By Luc Devroye[1] and Gábor Lugosi[2]

*McGill University and ICREA and Pompeu Fabra University*

It is shown that functions defined on $\{0, 1, \ldots, r - 1\}^n$ satisfying certain conditions of bounded differences that guarantee sub-Gaussian tail behavior also satisfy a much stronger "local" sub-Gaussian property. For self-bounding and configuration functions we derive analogous locally subexponential behavior. The key tool is Talagrand's [*Ann. Probab.* **22** (1994) 1576–1587] variance inequality for functions defined on the binary hypercube which we extend to functions of uniformly distributed random variables defined on $\{0, 1, \ldots, r - 1\}^n$ for $r \geq 2$.

**1. Introduction.** Concentration inequalities for functions of independent random variables establish upper bounds for the tail probabilities of such functions under general "smoothness" conditions; see, for example, Talagrand [30, 31, 32], Ledoux [19, 20], Boucheron, Lugosi, Massart [7, 8], McDiarmid [23], and so on. In this paper we take a closer look at the distribution of certain functions of independent random variables and show that the tail distribution exhibits a sub-Gaussian (or subexponential) behavior in a stronger "local" sense in many cases when concentration inequalities predict a sub-Gaussian (subexponential) tail.

First we consider real-valued functions defined on the binary hypercube $f : \{0, 1\}^n \to \mathbb{R}$. If $X = (X_1, \ldots, X_n)$ is uniformly distributed on the hypercube, we are interested in the distribution of the random variable $f(X)$.

Received May 2006; revised November 2006.
[1]Supported by NSERC Grant A3456 and FQRNT Grant 90-ER-0291.
[2]Supported by the Spanish Ministry of Science and Technology and FEDER, Grant BMF2003-03324 and by the PASCAL Network of Excellence under EC Grant 506778.

Our starting point is the following inequality, due to Talagrand [29]:

$$\mathrm{Var}(f) \le \frac{9}{10} \sum_{i=1}^{n} \frac{\mathbf{E}(f(X)-f(X^{(i)}))^2}{1+\log(\sqrt{\mathbf{E}(f(X)-f(X^{(i)}))^2}/(\mathbf{E}|f(X)-f(X^{(i)})|)},$$

(1.1)

where $X^{(i)} = (X_1, \ldots, 1-X_i, \ldots, X_n)$ is obtained by flipping the $i$th bit of $X$ and $\mathrm{Var}(f)$ denotes the variance of the random variable $f(X)$. The constants shown here follow from a simple proof by Benjamini, Kalai and Schramm [5].

Note that (apart from numerical constants) Talagrand's inequality improves upon the well-known Efron–Stein inequality (see Efron and Stein [11], Rhee and Talagrand [27], Steele [28]):

$$\mathrm{Var}(f) \le \tfrac{1}{2} \sum_{i=1}^{n} \mathbf{E}(f(X)-f(X^{(i)}))^2.$$

In Section 2 we show how to use Talagrand's inequality to prove "local" sub-Gaussian concentration inequalities. As a simple example, we show that if $f: \{0,1\}^n \to \mathbb{R}$ is such that there exists a constant $v$ such that $\sum_{i=1}^{n}(f(x) - f(x^{(i)}))_+^2 \le v$, then for all $k = 1, 2, 3, \ldots$,

$$a_{k+1} - a_k \le c\sqrt{v/k},$$

where $a_k$ denotes a $1 - 2^{-k}$ quantile of $f(X)$ and $c$ is a universal constant. The main argument is based on an observation of Benjamini, Kalai and Schramm [5] who show how Talagrand's inequality may be used to obtain exponential concentration inequalities. Even though Benjamini, Kalai and Schramm do not mention the possibility of deriving local concentration inequalities, it is their argument which is at the basis of our proofs. The purpose of this paper is to elaborate on this argument and to derive local concentration inequalities under different conditions. In Sections 3, 4 and 5 various variants and extensions are introduced. In Section 3 local concentration inequalities are shown under different conditions that are satisfied for numerous natural examples such as *configuration functions* introduced by Talagrand [30]—for *self-bounding functions*, see Boucheron, Lugosi and Massart [7], Maurer [22] and McDiarmid and Reed [24].

In Section 4, Talagrand's inequality is extended from the binary hypercube to functions defined on $\{0, 1, \ldots, r-1\}^n$ under the uniform distribution. The main technical tool here is a suitable hypercontractive inequality proved by Alon, Dinur, Friedgut and Sudakov [2]. This extension allows us to generalize the results of Sections 2 and 3 to functions defined on $\{0, 1, \ldots, r-1\}^n$.

In Section 5 we illustrate the use of the results of Section 4 by considering two classical, structurally similar, problems. We derive local concentration inequalities for the cost of the minimum weight spanning tree of a complete graph with random uniform weights on the edges and also for the assignment problem.

**2. Functions with locally sub-Gaussian behavior.** First we consider functions $f:\{0,1\}^n \to \mathbb{R}$ which satisfy the following properties: for all $x = (x_1, \ldots, x_n) \in \{0,1\}^n$,

$$\sum_{i=1}^{n}(f(x) - f(x^{(i)}))_+^2 \leq v, \tag{2.1}$$

where $v$ is a positive constant. [Here and throughout the paper, $a_+ = \max(a, 0)$ and $a_- = \max(-a, 0)$ denote the positive and negative parts of the real number $a$.] Clearly, if $f$ is 1-Lipschitz under the Hamming distance, then $v \leq n$, but there are many interesting examples in which $v$ is significantly smaller than $n$. It is well known (see Ledoux [19], or Boucheron, Lugosi and Massart [7]) that for such functions

$$\mathbf{P}\{f(X) \geq \mathbf{E}f(X) + t\} \leq e^{-t^2/4v}. \tag{2.2}$$

Our basic result (Theorem 2.1) shows that tail quantiles of the random variable $f(X)$ are not far apart. In this sense, it is a local tail bound. For any $\alpha \in (0,1)$, define the $\alpha$-quantile of $f$ by

$$Q_\alpha = \inf\{z : \mathbf{P}\{f(X) \leq z\} \geq \alpha\}.$$

In particular, we denote the median of $f(X)$ by $\mathbf{M}f = Q_{1/2}$.

THEOREM 2.1. *Assume $f$ satisfies (2.1) and let $B = \max_{x,i}|f(x) - f(x^{(i)})|$. Then for all $b > a \geq \mathbf{M}f$,*

$$b - a \leq \sqrt{\frac{(72/5)v\mathbf{P}\{f(X) \in (a, b+B)\}}{\mathbf{P}\{f(X) \geq b\}\log(e^2/(2\mathbf{P}\{f(X) \in (a, b+B)\}))}}$$

$$\leq \sqrt{\frac{(72/5)v\mathbf{P}\{f(X) > a\}}{\mathbf{P}\{f(X) \geq b\}\log(e^2/(2\mathbf{P}\{f(X) > a\}))}}.$$

PROOF. Define the function $g_{a,b}:\{0,1\}^n \to \mathbb{R}$ by

$$g_{a,b}(x) = \begin{cases} b, & \text{if } f(x) \geq b, \\ f(x), & \text{if } a < f(x) < b, \\ a, & \text{if } f(x) \leq a. \end{cases}$$

First observe that

$$\mathrm{Var}(g_{a,b}(X)) \geq \frac{\mathbf{P}\{f(X) \geq b\}}{4}(b-a)^2.$$

On the other hand, we may use Talagrand's inequality to obtain an upper bound for the variance of $g_{a,b}(X)$. To this end, observe that

$$\mathbf{E}|g_{a,b}(X) - g_{a,b}(X^{(i)})|$$
$$= 2\mathbf{E}(g_{a,b}(X) - g_{a,b}(X^{(i)}))_+$$
$$= 2\mathbf{E}[(g_{a,b}(X) - g_{a,b}(X^{(i)}))_+ \mathbf{1}_{f(X) \in (a, b+B)}]$$
$$\text{(by the definition of } g_{a,b} \text{ and } B)$$
$$\leq 2\sqrt{\mathbf{E}(g_{a,b}(X) - g_{a,b}(X^{(i)}))_+^2}\sqrt{\mathbf{P}\{f(X) \in (a, b+B)\}}$$
$$\text{(by Cauchy–Schwarz)}$$
$$= \sqrt{2\mathbf{E}(g_{a,b}(X) - g_{a,b}(X^{(i)}))^2}\sqrt{\mathbf{P}\{f(X) \in (a, b+B)\}}.$$

On the other hand,

$$\sum_{i=1}^n \mathbf{E}(g_{a,b}(X) - g_{a,b}(X^{(i)}))^2$$
$$= 2\sum_{i=1}^n \mathbf{E}(g_{a,b}(X) - g_{a,b}(X^{(i)}))_+^2$$
$$= 2\mathbf{E}\left[\mathbf{1}_{f(X) \in (a, b+B)} \sum_{i=1}^n (g_{a,b}(X) - g_{a,b}(X^{(i)}))_+^2\right]$$
$$\leq 2v\mathbf{P}\{f(X) \in (a, b+B)\},$$

where in the last step we used the fact that (2.1) implies that

$$\sum_{i=1}^n (g_{a,b}(X) - g_{a,b}(X^{(i)}))_+^2 \leq \sum_{i=1}^n (f(X) - f(X^{(i)}))_+^2 \leq v.$$

Combining the lower bound for the variance with the upper bound obtained by Talagrand's inequality yields the claim. □

To make Theorem 2.1 more transparent, we state a simple corollary for quantiles of $f(X)$. Using $\mathbf{P}\{f(X) > Q_{1-\gamma}\} \leq \gamma$ and $\mathbf{P}\{f(X) \geq Q_{1-\delta}\} \geq \delta$, Theorem 2.1 implies the following bound for the distance between any two quantiles in the upper tail:

THEOREM 2.2. *Assume $f$ satisfies (2.1). Then for all $\delta < \gamma \leq 1/2$,*

$$Q_{1-\delta} - Q_{1-\gamma} \leq \sqrt{\frac{(72/5)v\gamma}{\delta \log(e^2/(2\gamma))}}.$$

In particular, by choosing $\gamma = 2^{-k}$ and $\delta = 2^{-(k+1)}$ for some integer $k \geq 1$ and introducing

$$a_k = Q_{1-2^{-k}},$$

we get

(2.3) $$a_{k+1} - a_k \leq \frac{12}{\sqrt{5}}\sqrt{\frac{v}{(k-1)\log 2 + 2}} \leq 4\sqrt{\frac{v}{k}}.$$

Summing over $k = 1, 2, \ldots, m-1$ and using $\sum_{k=1}^{m-1}(k-1)^{-1/2} \leq \int_0^{m-1} x^{-1/2}\,dx = 2\sqrt{m-1}$, we obtain

$$a_m \leq a_1 + 8\sqrt{v(m-1)},$$

recovering (up to a constant factor) the sub-Gaussian concentration inequality (2.2) for $f$. However, Theorem 2.2 shows a sub-Gaussian behavior in a significantly stronger sense. If $f(X)$ was a normal random variable with variance $v$, then one would have $a_k \sim \sqrt{2vk\log 2}$ and $a_{k+1} - a_k \sim \sqrt{v\log 2/k}$. This (up to a constant factor) is precisely of the form of the upper bound (2.3) for a general function $f$ satisfying (2.1). Thus, the whole quantile sequence $\{a_k\}$ is a *contraction* of that of a normal random variable of variance a constant times $v$. (We say that a sequence $\{x_n\}$ is a contraction of another sequence $\{y_n\}$ if for every $n = 1, 2, \ldots$, $|x_{n+1} - x_n| \leq |y_{n+1} - y_n|$.)

REMARK (C). Even though we offer explicit numerical constants in the inequalities derived throughout the paper, no optimality of these values is claimed. In fact, quite often we sacrifice better constants for convenience in the notation or for simpler arguments.

EXAMPLE (C). One of the main examples of a function satisfying (2.1) is Talagrand's *convex distance* (Talagrand [30]) defined as follows. Let $A \subset \{0,1\}^n$ and define $f$ as

$$f(x) = \sup_{\alpha \in [0,\infty)^n : \|\alpha\|=1} \inf_{y \in A} \sum_{i\,:\,x_i \neq y_i} |\alpha_i|$$

where $x = (x_1, \ldots, x_n)$ and $y = (y_1, \ldots, y_n)$. Talagrand shows that for any set $A$ with $\mathbf{P}\{X \in A\} \geq 1/2$,

$$\mathbf{P}\{f(X) \geq t\} \leq 2e^{-t^2/4}.$$

(Note that Talagrand's result is true in any product space with product measure.) It is shown by Boucheron, Lugosi and Massart [8] that $f$ satisfies (2.1) with $v = 1$. This implies that for all $k = 1, 2, 3, \ldots$,

$$a_{k+1} - a_k \leq \frac{4}{\sqrt{k}}.$$

EXAMPLE (L). Let $f(X)$ denote the largest eigenvalue of the adjacency matrix of a random graph $G(m, 1/2)$ on $m$ vertices such that each edge appears with probability $1/2$. Thus, $n = \binom{m}{2}$ and $X_i = 1$ if and only if edge $i$ is present in the graph. Füredi and Komlós [14] show that $f(X)$ is asymptotically normally distributed with expectation $m/2$ and variance $1/2$. Alon, Krivelevich and Vu [3] show that $f(x)$ satisfies (2.1) with $v = 4$ (see also Maurer [22]) and conclude that $a_k \leq \mathbf{M}f(X) + \sqrt{32(k+2)\log 2}$. Theorem 2.2 implies the nonasymptotic local sub-Gaussian estimate

$$a_{k+1} - a_k \leq \frac{8}{\sqrt{k}}$$

for $k = 1, 2, \ldots$. Note that Alon, Krivelevich and Vu [3] also prove a concentration result for the $r$th largest inequality of the form $a_k \leq \mathbf{M}f(X) + Cr\sqrt{k}$. Their argument may be combined with ours to obtain an analogous local concentration inequality.

EXAMPLE (R). Another example is a *Rademacher average* of the form

$$f(x) = \sup_{\alpha \in A} \sum_{i=1}^{n} \alpha_i(x_i - 1/2),$$

where $A \subset \mathbb{R}^n$ is a set of vectors $\alpha$ with $\|\alpha\| \leq 1$. It is easy to see that condition (2.1) is satisfied with $v = 1$.

REMARK (A). We note here that Talagrand proved his inequality (1.1) in a more general setup in which the components $X_i$ of $X$ are i.i.d. Bernoulli($p$) random variables for some $p \in (0, 1)$. In this more general case Theorem 2.2 becomes

$$Q_{1-\delta} - Q_{1-\gamma} \leq \sqrt{\frac{Cv\gamma}{\delta \log(1/(2\gamma)) \log 1/(p(1-p))}}$$

for some constant $C$.

One obtains a corollary of a slightly different flavor by choosing, in Theorem 2.1, $a = k$ and $b = k + 1$ for some integer $k \geq \mathbf{M}f$; Theorem 2.1 implies the following local lower bound for the distribution of $f$:

COROLLARY 2.1. *Assume $f$ satisfies (2.1). Then for all $k \geq \mathbf{E}f + \sqrt{4v \log 2}$,*

$$\frac{q_k}{\sum_{i \geq k+1} q_i} + 1 \geq \frac{5}{288} \frac{(k - \mathbf{E}f)^2}{v^2} + \frac{5}{72v} \log \frac{e^2}{2}$$

*where $q_k = \mathbf{P}\{f(X) \in [k, k+1)\}$.*

PROOF. This follows immediately by noting that, on the one hand by Theorem 2.1, for $k \geq \mathbf{M}f$,

$$\sum_{i \geq k} q_i \leq (72/5)v \left( q_k + \sum_{k+1 \leq i \leq k+B+1} q_i \right) \left( \log \frac{e^2}{2(q_k + \sum_{k+1 \leq i \leq k+B+1} q_i)} \right)^{-1}$$

$$\leq (72/5)v \left( q_k + \sum_{i \geq k+1} q_i \right) \left( \log \frac{e^2}{2(q_k + \sum_{i \geq k+1} q_i)} \right)^{-1}$$

so that

$$q_k + \sum_{i \geq k+1} q_i \geq \frac{e^2}{2} \exp\left( -(72/5)v \left( \frac{q_k}{\sum_{i \geq k+1} q_i} + 1 \right) \right).$$

By the concentration inequality (2.2), for all $k \geq \mathbf{E}f$,

$$q_k + \sum_{i \geq k+1} q_i = \mathbf{P}\{f(X) \geq k\} \leq e^{-(k - \mathbf{E}f)^2/4v}.$$

Since $\mathbf{M}f \leq \mathbf{E}f + \sqrt{4v \log 2}$, combining the upper and lower bounds implies the corollary. $\square$

REMARK (*Monotonicity of the tail*). An obvious corollary is that $q_{k+1} \leq q_k$ whenever $k \geq \mathbf{E}f + (25/\sqrt{5})v$.

In some applications, even though (2.1) is not satisfied, the similar condition

(2.4) $$\sum_{i=1}^n (f(x) - f(x^{(i)}))_-^2 \leq v$$

holds. For such cases the next analog of Theorem 2.1 is true. The proof is omitted as it is a straightforward modification. In Section 5 we present some applications of this result.

THEOREM 2.3. *Assume $f$ satisfies (2.4) and let $B = \max_{x,i} |f(x) - f(x^{(i)})|$. Then for all $b > a \geq \mathbf{M}f$,*

$$b - a \leq \sqrt{\frac{(72/5)v \mathbf{P}\{f(X) \in (a - B, b)\}}{\mathbf{P}\{f(X) \geq b\}}} \left( \log \frac{e^2}{2\mathbf{P}\{f(X) \in (a - B, b)\}} \right)^{-1}.$$

*In particular, for all $\delta < \gamma \leq 1/2$, by taking $a = Q_{1-\gamma} + B$ and $b = Q_{1-\delta}$, we have*

$$Q_{1-\delta} - Q_{1-\gamma} \leq B + \sqrt{\frac{(72/5)v\gamma}{\delta \log(e^2/(2\gamma))}}.$$

**3. Configuration functions.** In this section we consider functions defined on the binary hypercube. Just as in Section 1, let $f:\{0,1\}^n \to \mathbb{R}$ and assume that $X$ is uniformly distributed over $\{0,1\}^n$.

Often, the sum of the squared changes appearing in condition (2.1) cannot be bounded by a constant but it can be related to the value of the function itself. Consider the following conditions:

$$|f(x) - f(x^{(i)})| \le B \quad \text{for all } x \text{ and } i \quad \text{and}$$

(3.1)

$$\sum_{i=1}^n (f(x) - f(x^{(i)}))_+^2 \le \phi(f(x)),$$

where $\phi$ is a fixed nonnegative nondecreasing function defined on the reals. In many applications, such as for configuration functions, one may take $\phi$ to be the identity and in some others it has the form $\phi(u) = au + b$ (see Talagrand [30], Boucheron, Lugosi and Massart [7, 8] and Devroye [9]). For example, it is shown by Boucheron, Lugosi and Massart [8] (for various extensions see also Maurer [22], McDiarmid and Reed [24]) that if (3.1) is satisfied with $\phi(u) = u$ and $B \le 1$, then

$$\mathbf{P}\{f(X) \ge \mathbf{E}f(X) + t\} \le e^{-t^2/(2\mathbf{E}f(X) + 2t/3)}.$$

Boucheron, Lugosi and Massart [8] offer concentration inequalities for the case when $\phi(u) = cu^\alpha$ for some $\alpha \in (0, 2)$.

A straightforward modification of the proof of Theorem 2.1 yields the following:

THEOREM 3.1. *Assume $f$ satisfies (3.1) and let $B = \max_{x,i} |f(x) - f(x^{(i)})|$. Then for all $b > a \ge \mathbf{M}f$,*

$$b - a \le \sqrt{\frac{(72/5)\phi(b + B)\mathbf{P}\{f(X) > a\}}{\mathbf{P}\{f(X) \ge b\} \log(e^2/(2\mathbf{P}\{f(X) > a\}))}}.$$

*Also, for all $\delta < \gamma \le 1/2$,*

$$Q_{1-\delta} - Q_{1-\gamma} \le \sqrt{\frac{(72/5)\phi(Q_{1-\delta} + B)\gamma}{\delta \log(e^2/(2\gamma))}}.$$

*In particular, recalling the notation $a_k = Q_{1-2^{-k}}$,*

$$a_{k+1} - a_k \le 4\sqrt{\frac{\phi(a_{k+1} + B)}{k}}.$$

EXAMPLE (*Self-bounding functions*). In many interesting applications, $\phi$ may be taken to be the identity function and $B = 1$. These functions have been called *self-bounding*; see Boucheron, Lugosi and Massart [7], Maurer [22], McDiarmid and Reed [24]. In general, if $\phi(u)$ is linear, then by the above-mentioned concentration inequality, for all $k \geq \mathbf{E}f(X)$, $a_k \leq ck$, and therefore

$$a_{k+1} - a_k \leq C$$

where $c, C$ are constants. Thus, in this case the quantile sequence $\{a_k\}$ is a contraction of that corresponding to an *exponentially* distributed random variable with parameter $O(1)$, in a similar sense that functions satisfying (2.1) had a quantile sequence contracting a Gaussian quantile sequence.

EXAMPLE (*Longest increasing subsequences*). Let now $x = (x_1, \ldots, x_n) \in \{0, 1, \ldots, r-1\}^n$ and define $f(x)$ to be the length of the longest increasing subsequence of $x_1, \ldots, x_n$, that is, the largest positive integer $m$ for which there exist $1 \leq i_1 < \cdots < i_m \leq n$ such that $x_{i_1} \leq x_{i_2} \leq \cdots \leq x_{i_m}$. Tracy and Widom [33] and Johansson [18] showed that if $X$ is uniformly distributed over $\{0, 1, \ldots, r-1\}^n$, then $(f(X) - n/r)/\sqrt{2n/r}$ converges, in distribution, to a random variable whose distribution depends on $r$ (see also Its, Tracy and Widom [16]). In the binary case (i.e., when $r = 2$), $f(x)$ is the longest subsequence of the form $000 \cdots 00111 \cdots 11$, and Theorem 3.1 may readily be used. It is immediate to see that $f$ satisfies (3.1) with $B = 1$ and $\phi(u) = u$ and therefore Theorem 3.1 implies a nonasymptotic local subexponential concentration inequality. [To see why (3.1) is satisfied, fix a maximal increasing subsequence in $x$ and observe that $(f(x) - f(x^{(i)}))_+ = 0$ whenever $x_i$ is not in this maximal sequence.] The same inequality holds when $f(x) = \log_2 N(x)$ where $N(x)$ is the number of all increasing subsequences of $x$. The fact that $\log_2 N(x)$ satisfies (3.1) with $B = 1$ and $\phi(u) = u$ was observed by Boucheron, Lugosi and Massart [7]. If $r > 2$, one may use the results of Section 4 below to obtain analogous bounds.

REMARK (*Concentration inequalities*). The recursion for the sequence $\{a_k\}$ given by Theorem 3.1 allows one to derive concentration inequalities for general functions $\phi$. We illustrate this for the example when $\phi(u) \leq cu^\alpha$ for some $c > 0$ and $\alpha \in [0, 2]$. Then Theorem 3.1 implies, after some work, that there exist constants $C, t_0$ such that for $t \geq t_0$,

$$\mathbf{P}\{f(X) \geq t\} \leq \begin{cases} Ce^{-t^{2-\alpha}/C}, & \text{if } 0 \leq \alpha < 2, \\ Ce^{-(\log t)^2/C}, & \text{if } \alpha = 2. \end{cases}$$

The case $\alpha < 2$ has already been dealt with by Boucheron, Lugosi and Massart [8], but the $\alpha = 2$ case seems to be new.

**4. Functions defined on the $r$-ary hypercube.** The purpose of this section is to extend the results of Theorems 2.1, 2.2 and 2.3 to functions $f$ defined on the $r$-ary cube $\{0, 1, \ldots, r-1\}^n$, equipped with the uniform distribution. In order to do this, we need to generalize Talagrand's variance inequality to this case. In particular, we prove the following:

THEOREM 4.1. *Let $r \geq 2$ be a positive integer and let $f : \{0, 1, \ldots, r-1\}^n \to \mathbb{R}$ be a real-valued function. Suppose $X = (X_1, \ldots, X_n)$ is uniformly distributed on $\{0, 1, \ldots, r-1\}^n$. For $1 \leq i \leq n$, $0 \leq j \leq r-1$ and for each $x = (x_1, \ldots, x_n)$, denote $x_{i,j} = (x_1, \ldots, x_{i-1}, x_i \oplus j, x_{i+1}, \ldots, x_n)$ where $\oplus$ stands for addition modulo $r$. Writing*

$$\Delta_i f(x) = f(x) - \frac{1}{r} \sum_{j=0}^{r-1} f(x_{i,j}),$$

*we have*

$$\mathrm{Var}(f) \leq 10(\log C_r) \sum_{i=1}^{n} \frac{\mathbf{E}(\Delta_i f(X))^2}{1 + \log(\sqrt{\mathbf{E}(\Delta_i f(X))^2}/\mathbf{E}|\Delta_i f(X)|)},$$

*where $C_r = (9/2)r^3$ is the constant of Lemma 4.1 below.*

As a consequence, Theorems 2.1, 2.2, 2.3 and 3.1 may now be extended to functions defined on $\{0, 1, \ldots, r-1\}^n$ with the only difference that in the conditions on $f$, $f(x) - f(x^{(i)})$ is replaced by $\Delta_i f(x)$ and the upper bounds in all four theorems are multiplied by $(10/3)\sqrt{\log C_r}$. For example, we will use the following result in Section 5:

COROLLARY 4.1. *Assume $f : \{0, 1, \ldots, r-1\}^n \to \mathbb{R}$ is such that there exists $v > 0$ such that*

$$\sum_{i=1}^{n} (\Delta_i)_-^2 \leq v$$

*and let $B = \max_{x,i} |\Delta_i f(x)|$. Then for all $k = 1, 2, 3, \ldots$,*

$$a_{k+1} - a_k \leq B + 14\sqrt{\log C_r}\sqrt{\frac{v}{k}}.$$

The proof of Theorem 4.1 is analogous to Talagrand's [29] original argument which was based on the Beckner–Bonami hypercontractive inequality (see Bonami [6] and Beckner [4]) of Fourier analysis on the binary hypercube. Here we use an extension of this inequality to functions defined on $\{0, 1, \ldots, r-1\}^n$ due to Alon, Dinur, Friedgut and Sudakov [2] which we recall below.

For any $S = (S_1, \ldots, S_n) \in \{0, 1, \ldots, r-1\}^n$, define the function

$$u_S(x) = \omega^{\langle S, x \rangle},$$

where $\omega = e^{2\pi i/r}$ and $\langle S, x \rangle = \sum_{i=1}^n S_i x_i \bmod r$. It is easy to see (see [2]) that the $u_S$ form an orthonormal basis of the space of complex-valued functions defined over $\{0, 1, \ldots, r-1\}^n$. To simplify notation, we will write

$$\int f = \frac{1}{r^n} \sum_{x \in \{0,1,\ldots,r-1\}^n} f(x) \quad \text{and} \quad \|f\|_q = \left(\int f^q\right)^{1/q}.$$

Denote by

$$\widehat{f}(S) = \int f \overline{u}_S$$

the Fourier coefficients of $f$ where $\overline{u}_S$ stands for the complex conjugate of $u_S$. A key ingredient of the proof is the following hypercontractive inequality:

LEMMA 4.1 (*Alon, Dinur, Friedgut and Sudakov* [2]). *For any $f : \{0, 1, \ldots, r-1\}^n \to \mathbb{R}$ and $k = 1, \ldots, n$,*

$$\left\| \sum_{S : |S| \leq k} \widehat{f}(S) u_S \right\|_4 \leq C_r^k \left( \sum_{S : |S| \leq k} \widehat{f}(S)^2 \right)^{1/2},$$

*where $C_r = (9/2) r^3$.*

PROOF OF THEOREM 4.1. Writing $f_{i,j}(x) = f(x_{i,j})$, it is easy to see that $\widehat{f}_{i,j}(S) = \widehat{f}(S) \omega^{jS_i}$. Thus,

$$\frac{1}{r} \sum_{j=0}^{r-1} \widehat{f}_{i,j}(S) = \begin{cases} \widehat{f}(S), & \text{if } S_i = 0, \\ 0, & \text{if } S_i \neq 0, \end{cases}$$

and therefore the Fourier coefficients of $\Delta_i f$ satisfy

$$\widehat{\Delta_i f}(S) = \begin{cases} 0, & \text{if } S_i = 0, \\ \widehat{f}(S), & \text{if } S_i \neq 0. \end{cases}$$

This and Parseval's identity imply that

$$\mathrm{Var}(f) = \|f\|_2^2 - \left(\int f\right)^2 = \sum_{S \neq \overline{0}} \widehat{f}(S)^2 = \sum_{i=1}^n \sum_{S \neq \overline{0}} \frac{\widehat{\Delta_i f}(S)^2}{|S|},$$

where $|S|$ denotes the number of nonzero components of $S$ and $\overline{0}$ is the all-zero vector.

Thus, in order to prove the theorem, it suffices to show that for any $f:\{0,1,\ldots,r-1\}^n \to \mathbb{R}$,

$$\sum_{S \neq \vec{0}} \frac{\widehat{f}(S)^2}{|S|} \leq 10 \log C_r \frac{\|f\|_2^2}{1+\log(\|f\|_2/\|f\|_1)},$$

which is what we do in the remaining part of the proof. Fix $k \leq n$ and observe that

$$\sum_{S:|S|=k} \widehat{f}(S)^2 = \int \left(\sum_{S:|S|=k} \widehat{f}(S) u_S\right) f$$

$$\leq \left\|\sum_{S:|S|=k} \widehat{f}(S) u_S\right\|_4 \cdot \|f\|_{4/3} \qquad \text{(by Hölder)}$$

$$\leq C_r^k \left(\sum_{S:|S|=k} \widehat{f}(S)^2\right)^{1/2} \cdot \|f\|_{4/3} \qquad \text{(by Lemma 4.1).}$$

This implies

$$\sum_{S:|S|=k} \widehat{f}(S)^2 \leq C_r^{2k} \|f\|_{4/3}^2$$

and we have, for all positive integers $m$,

$$\sum_{S:1 \leq |S| \leq m} \frac{\widehat{f}(S)^2}{|S|} \leq \|f\|_{4/3}^2 \sum_{k=1}^m \frac{C_r^{2k}}{k} \leq K \frac{C_r^{2m}}{m} \|f\|_{4/3}^2$$

where $K = \frac{36^2/2}{36^2/2 - 1}$. At the last step we used the fact that $C_r \geq 36$ and therefore $C_r^{2(k+1)}/(k+1) \geq (36^2/2) C_r^{2k}/k$. Now we may write

$$\sum_{S \neq \vec{0}} \frac{\widehat{f}(S)^2}{|S|} = \sum_{S:1 \leq |S| \leq m} \frac{\widehat{f}(S)^2}{|S|} + \sum_{S:|S|>m} \frac{\widehat{f}(S)^2}{|S|}$$

$$\leq K \frac{C_r^{2m}}{m} \|f\|_{4/3}^2 + \frac{1}{m+1} \sum_{S:|S|>m} \widehat{f}(S)^2$$

$$\leq \frac{1}{m+1} (2K C_r^{2m} \|f\|_{4/3}^2 + \|f\|_2^2).$$

Now we choose $m$ as the largest integer such that $C_r^{2m} \|f\|_{4/3}^2 \leq e^{1/3} \|f\|_2^2$ so that

$$m+1 \geq \frac{\log(e^{1/3} \|f\|_2/\|f\|_{4/3})}{\log C_r}$$

and

$$\sum_{S \neq \vec{0}} \frac{\widehat{f}(S)^2}{|S|} \leq \frac{1}{m+1} \cdot (2K+1)\|f\|_2^2 \leq \frac{(2K+1)\|f\|_2^2 \log C_r}{\log(e^{1/3}\|f\|_2/\|f\|_{4/3})}.$$

The proof is finished by observing that

$$\int f^{4/3} = \|f^{3/2}\|_{8/9}^{8/9} \leq \|f^{3/2}\|_1^{8/9} \leq (\|f\|_1 \cdot \|f\|_2^2)^{4/9}$$

by the Cauchy–Schwarz inequality, and therefore $\|f\|_{4/3}^3 \leq \|f\|_1 \cdot \|f\|_2^2$, which is equivalent to

$$\frac{e\|f\|_2}{\|f\|_1} \leq \left(\frac{e^{1/3}\|f\|_2}{\|f\|_{4/3}}\right)^3. \qquad \square$$

REMARK (*Logarithmic Sobolev inequalities*). An alternative route, yielding better numerical constants than Lemma 4.1, would be to use a sharp logarithmic Sobolev inequality of Diaconis and Saloff-Coste ([10], Theorem A.1) which implies hypercontractivity by Gross' theorem; see [15].

**5. Minimum weight spanning tree and the assignment problem.** In this section we derive local concentration bounds for two classical problems: the minimum weight spanning tree and the assignment problem. In these examples the random variables of interest are functions of independent random variables uniformly distributed in $[0,1]$. By simple discretization we may approximate them by functions defined over $\{0, 1, \ldots, r-1\}^n$ and use the result of the previous section. Since in Corollary 4.1 the dependence on $r$ is only logarithmic, we may take $r$ to be quite large (proportional to $n$ in these cases) and still obtain meaningful results.

Concentration inequalities for both cases may be derived, for example, by Talagrand's [30] results. In fact, Talagrand works out the case of the assignment problem. In order to conveniently use general concentration inequalities, Talagrand uses a truncation argument, a technique we also adopt below. Interestingly, the proofs in both examples below are identical and use simple structural properties of the function at hand.

EXAMPLE (*Minimum weight spanning tree*). Consider the random variable $T_m$ defined as the sum of weights on the minimum spanning tree of the complete graph $K_m$ with independent uniformly distributed (on $[0,1]$) weights $Y_{i,j}$ ($1 \leq i < j \leq m$) on the edges. A classical result of Frieze [13] shows that $\lim_{m \to \infty} \mathbf{E} T_m = \zeta(3)$. Janson [17] and Wästlund [36] prove that if the edge weights are exponentially distributed with parameter 1, then $\sqrt{m}(T_m - \zeta(3))$ converges, in distribution, to a centered normal random variable with variance $6\zeta(4) - 4\zeta(3)$. Here we study the related random variable $\overline{T}_m$ obtained when the $Y_{i,j}$ are replaced by $\min(Y_{i,j}, \delta_m)$ where $\delta_m > 0$ is a small positive number. Note that if $\delta_m = c\log m/m$ for some $c > 1$, then $T_m = \overline{T}_m$ with high probability. In order to see this just observe that $T_m \neq \overline{T}_m$ implies that the largest edge weight in the minimum spanning tree is greater than $\delta_m$. But this is just the probability that the random graph $G(m, \delta_m)$ is not connected which is at most $2(e^{m^{(1-c)/2}} - 1) + 2^{m+1}m^{-(c-1)m/4}$ (see Erdős and Rényi [12] and Palmer [26]), which is at most $4m^{-c/4}$, if $c \geq 2$.

To be able to use Corollary 4.1, we need to approximate $T_m$ by a function defined on $\{0, 1, \ldots, r-1\}^n$ under the uniform distribution where $n = \binom{m}{2}$. In order to do this, we replace the random variables $Y_{i,j}$ by their "discretized" approximation $\lfloor rY_{i,j} \rfloor / r$. If we denote the cost of the minimum spanning tree defined by the edge costs $\min(\lfloor rY_{i,j}\rfloor/r, \delta_m)$ by $\widetilde{T}_m$, then clearly $|\overline{T}_m - \widetilde{T}_m| \leq m/r$. The random variable $\widetilde{T}_m$ may now be considered as a function of $n = \binom{m}{2}$ independent variables $X_{i,j}$, all uniformly distributed on $\{0, 1, \ldots, r-1\}$, by defining $\lfloor rY_{i,j} \rfloor = X_{i,j}$. Now we may use Corollary 4.1. Clearly, we may take $B = \delta_m$. On the other hand,

$$\sum_{1 \leq i < j \leq m} (\Delta_{i,j})^2_- \leq m\delta_m^2$$

and therefore, denoting by $\widetilde{a}_k$ the $1 - 2^{-k}$-quantile of $\widetilde{T}_m$, we obtain

$$\widetilde{a}_{k+1} - \widetilde{a}_k \leq \delta_m + 14\sqrt{\frac{m\delta_m^2}{k}}\sqrt{\log(9r^3/2)}.$$

This, in turn, implies that if $\overline{a}_k$ denotes the $1 - 2^{-k}$-quantile of $\overline{T}_m$, then, for all $k = 1, 2, 3, \ldots$,

$$\overline{a}_{k+1} - \overline{a}_k \leq 2m/r + \delta_m + 14\sqrt{\frac{m\delta_m^2}{k}}\sqrt{\log(9r^3/2)}.$$

By choosing, say, $r = m^2$ and $\delta_m = c\log m/m$ for some constant $c > 1$, we obtain

$$\overline{a}_{k+1} - \overline{a}_k \leq C\left(\sqrt{\frac{\log^3 m}{mk}} + \frac{\log m}{m}\right)$$

for a constant $C$ depending on $c$ only. This inequality shows a local sub-Gaussian behavior whenever $k \leq m\log m$. It may be regarded as a local nonasymptotic version of the limit theorem of Janson and Wästlund, up to the logarithmic factors we needed to give up for technical reasons. For larger values of $k$ the second term dominates the first one, which corresponds to a subexponential behavior in the far tail. We do not know if this term is necessary. In order to convert this into a useful bound for the original problem $T_m$, one needs to choose $c$ so large that the bound $\mathbf{P}\{T_m \neq \overline{T_m}\} \leq 4m^{-c/4}$ does not dominate $2^{-k}$. Choosing $c = \max(2, 4(k+2)\log 2/\log m)$, one obtains

$a_{k+1} - a_{k-1}$
$$\leq \begin{cases} \dfrac{2}{m} + \dfrac{4(k+2)\log 2}{m} + 56\log 2\sqrt{\dfrac{3(k+2)}{m}\log \dfrac{9m^6}{2}}, \\ \hspace{6cm} \text{if } k+2 > \dfrac{\log m}{2\log 2}, \\ \dfrac{2}{m} + \dfrac{2\log m}{m} + 28\sqrt{\dfrac{\log^2 m}{km}\log \dfrac{9m^6}{2}}, \hspace{1cm} \text{if } k+2 \leq \dfrac{\log m}{2\log 2}. \end{cases}$$

In order to compare this local bound to concentration inequalities, note that Theorem 7 of Boucheron, Lugosi and Massart [8] implies that $\mathbf{P}\{\overline{T}_m \geq \mathbf{E}\overline{T}_m + t\} \leq e^{-t^2 m/(4(e-1)c^2 \log^2 m)}$, or in other words, that

$$\overline{a}_k \leq \mathbf{E}\overline{T}_m + \sqrt{\dfrac{kc^2 \log^2 m \log 2}{4(e-1)m}}.$$

Again, choosing $c = \max(2, 4(k+2)\log 2/\log m)$, one obtains

$$a_{k-1} \leq \mathbf{E}\overline{T}_m + \sqrt{\dfrac{\log 2}{e-1}} \max\left(\sqrt{\dfrac{k\log^2 m}{m}}, \sqrt{\dfrac{4(k+2)^3}{m}}\right).$$

By summing the "local" inequality in $k$, one obtains a concentration inequality that is only slightly weaker than the one derived here, as we get an extra factor of $\sqrt{\log m}$. This is due to the approximation by discretization, necessary to apply Corollary 4.1.

EXAMPLE (The assignment problem). In the assignment problem, given an $m \times m$ array $\{Y_{i,j}\}_{m \times m}$ of independent random variables distributed uniformly on $[0,1]$, one considers the random quantity

$$Z_m = \min_\pi \sum_{i=1}^m Y_{i,\pi(i)}$$

where the minimum is taken over all permutations $\pi$ of $\{1, \ldots, m\}$. Culminating a long series of partial results, Aldous [1] shows that $\lim_{m \to \infty} \mathbf{E} Z_m = \zeta(2)$. In the case when the $Y_{i,j}$ are exponentially distributed with parameter 1, Linusson and Wästlund [21] and Nair, Prabhakar and Sharma [25] independently prove that for all $m$, $\mathbf{E} Z_m = \sum_{i=1}^m i^{-2}$. See also Wästlund [34]. Wästlund [35] also derives an explicit formula for the variance of $Z_m$. In particular, he proves that $\mathrm{Var}(Z_m) = 4(\zeta(2) - \zeta(3))/m + O(m^{-1/2})$. Talagrand [30] proves (in the uniform model) an exponential concentration inequality very similar to the one described for the minimum weight spanning tree above.

In fact, in order to get local concentration inequalities, we may proceed exactly as we did in the previous example: first we replace the $Y_{i,j}$ by the truncated variables $\min(Y_{i,j}, \delta_m)$. If $\overline{Z}_m$ denotes the cost of the optimal assignment based on the truncated variables, then Proposition 10.3 of Talagrand [30] implies that there exists a constant $K$ such that $\mathbf{P}\{Z_m \neq \overline{Z}_m\} \leq e^{-m\delta_m/K}$, an inequality that is completely analogous to the one we used in the study of the minimum weight spanning tree. Second, we use the discretized approximation of the truncated variables. Then just as for the minimum weight spanning tree, we may take $B = \delta_m$ in Corollary 4.1 and observe that

$$\sum_{1 \leq i < j \leq m} (\Delta_{i,j})_-^2 \leq m\delta_m^2,$$

which leads to inequalities completely analogous to those obtained for the minimum weight spanning tree example above. In particular, if $a_k$ denotes the $1 - 2^{-k}$ quantile of $Z_m$, then there exists a constant $C$ such that

$$a_{k+1} - a_{k-1} \leq C \max\biggl(\frac{k}{m} + \sqrt{\frac{k \log m}{m}}, \frac{\log m}{m} + \sqrt{\frac{\log^3 m}{km}}\biggr).$$

SCHOOL OF COMPUTER SCIENCE
MCGILL UNIVERSITY
MONTREAL, QUEBEC
CANADA H3A 2A7
E-MAIL: luc@cs.mcgill.ca

ICREA AND DEPARTMENT OF ECONOMICS
POMPEU FABRA UNIVERSITY
RAMON TRIAS FARGAS 25-27
08005 BARCELONA
SPAIN
E-MAIL: lugosi@upf.es